\documentclass[11pt]{article}

\usepackage{amsfonts,amsmath,latexsym}
\usepackage{amstext,amssymb}
\usepackage{amsthm}
\usepackage{color}

\newcommand {\R}{\mathbb{R}}

\newcommand {\C}{\mathbb{C}}
\newcommand {\D}{\mathbb{D}}

\newcommand {\calC}{\mathcal{C}}
\newcommand {\calCp}{\mathcal{C}_p}
\newcommand {\calS}{\mathcal{S}}

\newtheorem {thm} {Theorem}

\newcommand{\beq}{\begin{equation}}
\newcommand{\eeq}{\end{equation}}
\numberwithin{equation}{section}

\def\labl#1{\label{#1}}
\begin {document}

\title{ 
{\bf Two Isoperimetric Inequalities for  
the Sobolev Constant}\footnote{AMS subject classification: 35J}}
\author{Tom Carroll 
 and Jesse Ratzkin  \\
University College Cork and University of Cape Town\\
{\tt t.carroll@ucc.ie} 
and {\tt jesse.ratzkin@uct.ac.za}} 
\maketitle

\begin {abstract} 
\noindent In this note we prove two isoperimetric inequalities for the sharp 
constant in the Sobolev embedding and its associated extremal function. 
The first such inequality is a variation on the classical Schwarz Lemma from 
complex analysis, similar to recent inequalities of Burckel, Marshall, Minda, 
Poggi-Corradini, and Ransford, while the second generalises an 
isoperimetric inequality for the first eigenfunction of the Laplacian 
due to Payne and Rayner. 
\end {abstract}

\section{Introduction and statement of results}

Let $p \geq 1$ and let $D$ be a bounded domain in the plane, equipped with the 
usual Lebesque area measure $dA$. The 
well-known Sobolev embedding theorem  states
that there exists a constant $k = k(p,D)$ such that 
\begin {equation} \label {sobolev-emb} 
u \in W^{1,2}_0(D) \Rightarrow  \|u\|_{L^p(D)} \leq k\|\nabla u \|_{L^2(D)},  
\end {equation} 
so that $W^{1,2}_0 (D) \hookrightarrow L^p(D)$.  
(This embedding is always compact for finite $p$ in two dimensions.) 
The sharp constant in the Sobolev embedding, 
\begin {equation} \label{sobolev-const}
\calS_p(D) = \sup \left \{ \frac {\left ( \int_D u^p\, dA\right )^{1/p}}
{\left ( \int_D |\nabla 
u|^2\, dA \right )^{1/2}} : u \in W^{1,2}_0(D), u \not \equiv 0 
\right \},\end {equation} 
is an important and well-studied constant related to the geometry and function 
theory of $D$. 
For example, van den Berg \cite{vdB} has recently derived asymptotic expansions for $\calS_p(D)$
as $p$ tends to $\infty$.
For our purposes, it is slightly more convenient to 
work with 
\begin {equation} \label {Cp}
\calC_p(D) = (\calS_p(D))^{-2} = \inf \left \{ 
\frac{\int_D |\nabla u|^2\, dA}{\left (
\int_D u^p\, dA\right )^{2/p} } : u \in W^{1,2}_0(D), u\geq 0, u \not 
\equiv 0 \right \}.
\end {equation}

The sharp Sobolev constant $\calS_p$ and its associated extremal 
functions are the subject of a vast literature which we make 
no attempt to survey here. Suffice to say that they are still 
the subject of much intense research and, in particular, there is 
a well-established link between \eqref{sobolev-emb} and the 
isoperimetric inequality (see \cite{DH} and references therein). 

Sobolev constants appear in a different guise in the theory of elasticity and 
the theory of vibrating membranes.
The constant  $\calC_2(D)$ minimises the Rayleigh-Ritz quotient and 
corresponds to  the first Dirichlet eigenvalue  $\lambda(D)$ 
of the Laplacian for $D$ (or, more generally, the bottom of the spectrum). 
The constant $\calC_1(D)$ is related to the torsional
rigidity $P(D)$ by $\calC_1(D) = 4/P(D)$. 
In this sense, $\calCp(D)$ interpolates between the torsional rigidity of a domain 
and its principal frequency as $p$ ranges from 1 to 2.

By the compactness of the embedding $W^{1,2}_0(D) \hookrightarrow L^p(D)$, a minimiser $\phi$ 
of the ratio defining $\calC_p(D)$ exists. 
The Euler-Lagrange equation tells us that $\phi$ solves the boundary value problem 
\beq\labl{p1.1}
\Delta \phi + \Lambda\, \phi^{p-1} = 0, \qquad 
\left. \phi \right|_{\partial D} = 0, 
\eeq
and an integration by parts agrument (\cite[Lemma~2]{CR}) 
shows that $\Lambda$ and $\calC_p(D)$ are related by 
\beq\labl{1.2}
\calCp(D) = \Lambda\left(\int_D \phi^p\,dA\right)^{(p-2)/p}.
\eeq

Our aim in this paper is two-fold: to prove a reverse H\"older inequality for 
the minimiser in \eqref{Cp} which generalises the inequality of Payne 
and Rayner \cite{PR} for the first Dirichlet eigenfunction of the Laplacian 
and, second, to 
prove a Schwarz Lemma for the constant $\calCp(D)$. 
In the rest of this introduction, we describe each of these inequalities 
in detail, and discuss how each is an isoperimetric inequality.

\subsection{A variation on the Schwarz Lemma} 
Let $f$ be a complex analytic function in the unit disc $\D = \{z: 
\vert z \vert < 1\}$ 
in the complex plane $\C$ with $f(0)=0$. 
Let $M(r)  = \max \{ \vert f(z)\vert \colon \vert z \vert =r\}$, 
$0< r < 1$,  be the maximum modulus of $f$ on the circle of 
radius $r$. Then $\lim_{r \to 0+} M(r)/r = \vert 
f'(0)\vert$ and $\lim_{r \to 1-} M(r)/r = \Vert f\Vert_\infty$. 
The classical Schwarz Lemma states that $M(r)/r$ increases with  $r$  
and that, if $M(r)/r$ takes the same value at two distinct $r \in [0,1]$, 
then $f$ is linear. 

Burckel, Marshall, Minda, Poggi-Corradini, and 
Ransford \cite{BMMPR} recently 
proved versions of the Schwarz Lemma for diameter, 
logarithmic capacity, and area. 
Laugesen and Morpurgo \cite{LM} proved a version of Schwarz's Lemma 
for principal frequency, 
in fact for the sum of the values of a convex increasing function applied 
to the 
reciprocals of the first $n$ Dirichlet eigenvalues of the Laplacian.
More recently, Betsakos has obtained counterparts of the results of 
Burckel {\it et al.}\ for quasiregular mappings \cite{Betsakos}.
We adapt the method of Laugesen and Morpurgo to 
prove a version of Schwarz's lemma for $\calCp$. 
\begin {thm} \labl{thm1}
Let $f$ be a conformal mapping of the unit disk $\D$ and let $p\geq 1$. 
The function
\beq\labl{1.1}
r \mapsto \frac{\calC_p\big( f(r\D) \big)}{\calC_p(r\D)}
 	= \frac{r^{4/p}}{\calC_p(\D)} \calC_p\big( f(r\D) \big), 
\quad 0 < r < 1,
\eeq
is strictly decreasing unless $f$ is linear (in which case this function 
is constant). 
Moreover, if $p\leq 2$, the reciprocal of this function is a convex 
function of $\log r$.
\end {thm}
The simple scaling law, $\calC_p(rD) = r^{-4/p}\calC_p(D)$, was used in 
\eqref{1.1}.
Using monotonicity and taking the limit of the right hand side of \eqref{1.1} 
as $r\rightarrow 0^+$, we obtain
\[
\calCp\big(f(\D)\big) \leq \calCp(\D) \, \vert f'(0)\vert^{-4/p},
\] 
which generalises the eigenvalue estimate in Section 5.8 of P\'olya and 
Szeg\H o \cite{PS}. Using monotonicity, we also see that the limit as
$r \rightarrow 1^-$ of the right hand side of \eqref{1.1} 
exists, though it might be zero, as it is in the case $f(z) = (1-z)/(1+z)$. 
The case $p=2$ is a special case of the result of Laugesen and Morpurgo 
\cite[Section~11]{LM} mentioned above.
The case $p=1$ is a Schwarz Lemma for torsional rigidity, which answers 
a question left open in \cite{BMMPR}.

One can interpret these variations on the Schwarz Lemma as 
dynamic isoperimetric inequalities. As $r$ increases, geometric 
quantities such as diameter and area of $f(r\D)$ increase as well, 
and the results of \cite{BMMPR} state they increase at least as rapidly 
as the case when $f(r\D)$ is a disk. Moreover, 
the case of equality only occurs if $f(r\D)$ is a disk with centre $f(0)$. Similarly, 
Theorem \ref{thm1} states that as $r$ increases the quantity 
$\calC_p(f(r\D))$ decreases at least as quickly as that of the 
disk $r\D$, with equality if and only if $f(r\D)$ is a disk with centre $f(0)$. 

\subsection{A reverse H\"older inequality for the minimiser of $\calC_p$}
Payne and Rayner \cite{PR} published a reverse H\"older inequality for the 
first Dirichlet eigenfunction $\phi$ of the Laplacian -- the minimiser of 
$\calC_2(D)$ -- that they had discovered some years before. 
Since they \lq saw at that time little use for the inequality\rq\  
and hoped to establish analogous results in higher dimensions 
they did not publish their result until 1972 coinciding with the publication 
of work by Sperb \cite{Sperb} in which he made 
\lq extensive use of this inequality\rq. 
The Payne-Rayner inequality for the first Dirichlet eigenfunction $\phi$ is
\beq\labl{PR}
\left( \int_D \phi\,dA\right)^2 \geq \frac{4\pi}{\lambda(D)} 
\int_D \phi^2\,dA.
\eeq  
Several extensions of the original equality have since been obtained, including
extensions to higher dimensions by Payne and Rayner \cite{PR2} and 
Kohler-Jobin \cite{KJ}, to equations more general than 
$\Delta u + \lambda u = 0$ with Dirichlet boundary conditions as in 
Alvino, Ferrone and Trombetti \cite{Alvino}, Chiti \cite{Chiti} and 
Mossino \cite{Moss}  and, most recently, to the  setting of minimal surfaces by 
Wang and Xia \cite{Wang}.  We prove the following extension. 
\begin{thm} \labl{thm2}Let $n \geq 3$ and let $\Sigma \subset \R^n$ be 
a compact, embedded, minimal surface with Lipschitz, weakly connected 
boundary $\partial \Sigma$. For $p \geq 1$ let  $\calCp(\Sigma)$ be given 
by \eqref{Cp} and let 
$\phi$ be a minimiser for $\calCp(\Sigma)$. Then 
\beq\labl{pPR}
\left( \int_\Sigma \phi^{p-1}\,dA\right)^2 \geq \frac{8\pi}{p\,\calCp(\Sigma)} 
\left(\int_\Sigma \phi^p\,dA\right)^{2-2/p}
\eeq  
Equality holds in \eqref{pPR} if and only of $\Sigma$ is a flat disk in an 
affine plane.
\end{thm}
We explain the relevant terminology for minimal surfaces, including the 
definition of a weakly connected boundary, at the beginning of 
Section~\ref{sec3}. Also, observe that Theorem~\ref{thm2} 
includes the setting when $\Sigma$ is a bounded planar domain as a 
special case. 

The line of thought behind our work in \cite{CR} was to obtain results for 
$\calCp$ that would interpolate between the corresponding results for 
torsional rigidity when $p=1$ and principal frequency when $p=2$, 
with the goal of both unifying and generalising such results. 
It is instructive from this point of view to note that, 
in the case of a bounded planar domain, the inequality
\eqref{pPR} reduces, of course, to the classical Payne-Rayner inequality 
when $p=2$, 
and becomes the Saint Venant inequality $2\pi P(D) \leq A(D)^2$, 
where $A(D)$ is the area of $D$, when $p=1$. 
The Saint Venant inequality is the isoperimetric inequality for torsional 
rigidity and states that, among all domains in the plane of given area, a 
disk of that area has the largest torsional rigidity. It was first proved 
by P\'olya in 1948.

We can 
rewrite \eqref{pPR} more geometrically by giving $\Sigma$ the 
(singular) conformal metric $\widetilde{ds} = 
\vert \nabla \phi \vert\, ds$, where 
$ds$ is the arc-length element on $\Sigma$ making its inclusion in 
Euclidean space an isometric embedding. 
The length $L$ of $\partial \Sigma$ in this metric is 
\[
L = \int_{\partial \Sigma} \vert \nabla \phi \vert \,ds
= -\int_{\partial \Sigma} \frac{\partial \phi}{\partial \eta}\,ds
=- \int_\Sigma\Delta \phi\, dA 
= \Lambda \int_\Sigma \phi^{p-1}\,dA,
\]
where $\partial/ \partial \eta$ denotes the outward normal derivative and 
we have used \eqref{p1.1}. By \eqref{1.2}, 
\[
L = \frac{\calCp(\Sigma)}{\left( \int_\Sigma\phi^p\,dA\right) ^{(p-2)/p}}
\int_\Sigma \phi^{p-1}\,dA.
\]
The area $A$ of the surface $\Sigma$ in this conformal metric is, using the 
fact that $\phi$ is a minimiser for $\calCp(\Sigma)$, 
\[
A = \int_\Sigma \vert \nabla \phi \vert^2\,dA = \calCp(\Sigma) 
\left(\int_\Sigma \phi^p\,dA\right)^{2/p}.
\]
In terms of the length $L$ of the boundary and the area $A$ of the surface 
with respect to the conformal metric induced by the length of the gradient 
of the minimiser $\phi$ of $\calCp(\Sigma)$, the inequality \eqref{pPR} 
becomes 
\beq\labl{isoPR}
L^2 \geq \frac{8\pi}{p}\,A,
\eeq
with equality if and only if $\Sigma$ is a flat disk in an affine plane. 
As previously pointed out by the authors \cite{CRSchwarz}, this inequality 
has, in the case $p=2$ 
of the first Dirichlet eigenfunction for the Laplacian, the same form as 
the classical isoperimetric inequality. 

\subsection {Organisation of the remainder of this paper} 
We prove Theorem~\ref{thm1} in Section~\ref{sec2} by a suitable adaptation of 
the method of Laugesen and Morpurgo \cite[Section~11]{LM}. Next we 
prove Theorem~\ref{thm2} in Section~\ref{sec3}, the proof being modelled on 
Payne and Rayner's original proof for the case $p=2$. 

\bigskip\noindent{\sc Acknowledgements:} We first learned about variations 
on Schwarz's Lem\-ma from Pietro Poggi-Corradini during a 
\textsl{Summer School in Conformal Geometry, Potential Theory, and Applications\/} 
held at NUI Maynooth in June 2009. We would like to thank Pietro 
for interesting discussions on the subject and the organizers of the conference for 
providing a stimulating venue for these discussions. 
Michiel van den Berg kindly told us about the Payne-Rayner inequality 
when we needed precisely this result in connection with a different proof 
(see \cite{CRSchwarz}) of a Schwarz Lemma for principal frequency.  T.\,C. is
partially supported by the ESF as part of the \lq Harmonic and Complex 
Analysis and Applications\rq\ programme, and J.\,R. is partially supported by the University 
of Cape Town Research Committee.  

\section {Proof of a Schwarz Lemma for $\calCp$}
\labl {sec2} 
For $\zeta\in\D$, let $g_\zeta$ be the conformal map of the disk $\D$ 
defined by $g_\zeta(z) = f(\zeta z)/\zeta$. Set $\Omega_\zeta = g_\zeta(\D)$. 
Note that $\Omega_\zeta$ is a rotation of $\Omega_{\vert \zeta \vert}$ so that 
$\calCp(\Omega_\zeta)$ depends only on $\vert  \zeta \vert$. 
Moreover,  $\Omega_r = f( r\D )/r$ for $0< r < 1$ so that, by the scaling law, 
$\calCp(\Omega_r) = r^{4/p}\calCp\big( f( r\D) \big)$. 

Set
\[
\Psi(\zeta) = \frac{1}{\calCp(\Omega_\zeta)}, \quad \zeta \in \D. 
\]
Following the method of Laugesen and Morpurgo, we show that 
$\Psi^{p/2}$ is subharmonic in $\D$. 
The integral means of this function are therefore both increasing and 
log-convex (see \cite[Theorem 2.12]{HK}). 
Since $\Psi^{p/2}$ is a radial function, its integral mean over the circle 
centre 0 and radius $r$ equates to $\Psi^{p/2}(r)$, for $0<r<1$, which 
is therefore both increasing and log-convex. 
It follows that $\Psi(r)$ is increasing, so that $\calCp(\Omega_r)$ is decreasing.
If, in addition, $p\leq 2$ then $2/p\geq 1$ so that 
the function $\Psi(r)$ is log-convex on $0< r <1$. 

The proof will therefore be complete once we show that $\Psi^{p/2}$ is 
subharmonic in the disk. First,
\[
\Psi^{p/2}(\zeta) = \sup  \left\{ 
\int_{\Omega_\zeta} u^p dA: 
u \in L^p(\Omega_\zeta) \cap W^{1,2}_0(\Omega_\zeta), 
u \geq 0, 
\int_{\Omega_\zeta}\vert \nabla u|^2 dA =1  \right\}.
\]
Write $\phi_\zeta$ for the extremal function on $\Omega_\zeta$ 
and $\psi_\zeta = \phi_\zeta \circ g_\zeta$ on $\D$. 
Since $g_\zeta'(z) = f'(\zeta z)$, $z \in \D$,  
\beq\labl{p2.1}
\Psi^{p/2}(\zeta) 
= \int_{\Omega_\zeta} \phi_\zeta^p\, dA
= \int_\D \psi_\zeta^p(z)\, \vert f'(\zeta z) \vert^2\,dA(z).
\eeq
Also, $1=\int_{\Omega_\zeta} \vert \nabla \phi_\zeta \vert^2\,dA = 
\int_\D \vert \nabla \psi_\zeta \vert^2\,dA$. 
In fact, if $\psi$ is any function on the disk for which 
$\int_\D \vert \nabla \psi \vert^2\,dA = 1$ and we set 
$\phi = \psi \circ g_\zeta^{-1}$ on $\Omega_\zeta$ then
\beq\labl{p2.2}
\int_{\Omega_\zeta} \vert \nabla \phi \vert^2\,dA = 1
\ \mbox{ and }\  
\int_{\Omega_\zeta} \phi^p\, dA
= \int_\D \psi^p(z)\, \vert f'(\zeta z) \vert^2\,dA(z).
\eeq
It then follows from \eqref{p2.1} and \eqref{p2.2} that
\beq\labl{p2.3}
\int_\D \psi^p(z)\, \vert f'(\zeta z) \vert^2\,dA(z)
\leq \int_\D \psi_\zeta^p(z)\, \vert f'(\zeta z) \vert^2\,dA(z)
\eeq
Fix $\zeta \in \D$ and let $0 < \rho < 1-\vert \zeta\vert$. 
First using \eqref{p2.1}, then using \eqref{p2.3} for the first inequality and the 
subharmonicity of $\vert f'\vert^2$ for the second, 
\begin{align*}
\frac{1}{2\pi} \int_0^{2\pi} \Psi^{p/2}(\zeta + &\rho e^{i\theta}) \,d\theta \\
& = \frac{1}{2\pi} \int_0^{2\pi}  \left( \int_\D \psi_{\zeta +\rho e^{i\theta}}^p(z)\, 
	\vert f'\big((\zeta + \rho e^{i\theta}) z\big) \vert^2\,dA(z) \right) \,d\theta\\
& \geq \frac{1}{2\pi} \int_0^{2\pi}  \left( \int_\D \psi_\zeta^p(z)\, 
	\vert f'\big((\zeta + \rho e^{i\theta}) z\big) \vert^2\,dA(z) \right) \,d\theta\\
& =  \int_\D \psi_\zeta^p(z)\, \left( \frac{1}{2\pi} \int_0^{2\pi} 
	\vert f'\big((\zeta + \rho e^{i\theta}) z\big) \vert^2\, d\theta\right) \,dA(z)\\
& \geq  \int_\D \psi_\zeta^p(z)\,  \vert f'(\zeta z) \vert^2 \,dA(z)\\
& =  \Psi^{p/2}(\zeta),
\end{align*}
thereby establishing the sub-mean value property. 
The last inequality is strict unless $f$ is linear, making $\Psi^{p/2}$ strictly subharmonic 
if $f$ is not linear. As shown by Laugesen and Morpurgo \cite[Page 104]{LM}, 
if a function is strictly subharmonic in the disk then its integral means are strictly increasing 
and strictly log-convex. This establishes the equality statement of Theorem~\ref{thm1}.

\section{Proof of a Payne-Rayner inequality}
\labl {sec3} 

In this section we prove a generalised Payne-Rayner inequality on a 
minimal surface. Some preliminary comments are in order. We consider a 
compact, two-dimensional, minimal surface $\Sigma \hookrightarrow \R^n$, 
with the induced metric, where $\Sigma$ has a Lipschitz boundary $\partial 
\Sigma$. The fact that $\Sigma$ is minimal means it is a critical point of 
the area function for variations which fix $\partial \Sigma$, or, equivalently, 
that the restriction of the coordinate functions to $\Sigma$ are all harmonic. 
Using the classical tools of the Riemannian geometry of surfaces, we define 
the usual gradient, 
divergence, and Laplace-Beltrami operators on $\Sigma$, and also 
$$\calC_p(\Sigma) = \inf \left \{ \frac{\int_\Sigma |\nabla u |^2 dA}{\left (
\int_\Sigma u^p dA \right )^{2/p} }: u \in W^{1,2}_0(\Sigma), u \not \equiv 0
\right \}$$
as before. The Sobolev embedding is still 
compact for $p \geq 1$, and so a minimiser $\phi$ exists and solves the 
boundary value problem 
$$\Delta_\Sigma \phi + \Lambda \phi^{p-1} = 0, \qquad \left. \phi 
\right |_{\partial \Sigma} =0,$$
where $\Lambda$ and $\calC_p(\Sigma)$ are related as in \eqref{1.2}. 

Li, Schoen, and Yau defined the notion of a weakly connected boundary 
$\partial \Sigma$ for a surface $\Sigma \subset \R^n$ in \cite{LSY}. The 
boundary $\partial \Sigma$ is weakly connected if there exists a rectangular 
coordinate system $\{x_1, \dots, x_n\}$ for $\R^n$ such that no coordinate 
hyperplane $\{ x_j = \textrm{constant} \}$ separates $\partial \Sigma$. They
prove that if $\Sigma$ is a compact minimal surface with Lipschitz, 
weakly connected boundary $\partial \Sigma$ then it satisfies the 
isoperimetric inequality $(L(\partial \Sigma))^2 \geq 4\pi A(\Sigma)$, 
with equality if and only if $\Sigma$ is a flat disk. 

Let $\{x_1, \dots, x_n\}$ be the usual rectangular coordinates on $\R^n$ 
and define the function 
$$f = \frac{1}{2} (x_1^2 + \cdots + x_n^2) $$
restricted to $\Sigma$. As a result of $\Sigma$ being minimal we have 
$\Delta_\Sigma f = 2$. 

\begin {proof} [Proof of Theorem \ref{thm2}]
Let $\phi_M = \max \{ \phi(p): p \in \Sigma\}$ and, for 
$0 \leq t \leq \phi_M$, define 
$$\Sigma(t) = \{ p \in \Sigma : \phi(p) \geq t \}, \qquad 
S(t) = \{p \in \Sigma: \phi(p) = t\}.$$
By Sard's theorem, $S(t) = \partial \Sigma(t)$ for almost 
every value of $t$. 

Now define 
$$H_0(t) = \int_{\Sigma(t)} \phi^{p-1}\, dA, \quad 
H_1 (t) = -\frac{p}{2} \int_{\Sigma(t)} \phi^{p-1} \langle 
\nabla \phi, \nabla f \rangle\, dA.$$
It will be useful to rewrite $H_1(t)$ as 
\begin {eqnarray} \label{min-H1} 
H_1(t) & = & -\frac{p}{2} \int_{\Sigma(t)} \phi^{p-1} \langle 
\nabla \phi, \nabla f \rangle\, dA = -\frac{1}{2} \int_{\Sigma(t)}
\langle \nabla (\phi^p), \nabla f\rangle\, dA  \nonumber \\
& = & \frac{1}{2} \int_{\Sigma(t)} \phi^p\, \Delta_\Sigma f\, dA - \frac{1}{2} 
\int_{\partial \Sigma(t)} \phi^p\, \frac{\partial f}{\partial \eta} 
\,ds  \nonumber \\ 
& = & \int_{\Sigma(t)} \phi^p\, dA - \frac{1}{2} \int_{\partial \Sigma(t)}
\phi^p\, \frac{\partial f}{\partial \eta}\, ds. 
\end {eqnarray} 

By the coarea formula, we have 
\begin {equation}\label{min0} 
H_0'(t) = \frac{d}{dt} \left(\int_t^{\phi_M} \tau^{p-1} \int_{S(\tau)}
\frac{ds}{|\nabla \phi|}\, d\tau\right) 
= -t^{p-1} \int_{S(t)} \frac{ds}{|\nabla \phi|}. 
\end {equation} 

Observe that $\eta = - \frac{\nabla \phi}{|\nabla \phi|}$ is the 
outward unit normal to $\Sigma(t)$. Hence, by the divergence theorem 
\begin {eqnarray} \label{min1} 
\int_{S(t)} |\nabla \phi| ds & = & - \int_{S(t)}\frac{1}{|\nabla \phi|}\,
 \langle \nabla \phi, -\nabla \phi \rangle\, ds = - \int_{\Sigma(t)} 
\Delta_\Sigma (\phi)\, dA \nonumber \\ 
& = & \Lambda \int_{\Sigma(t)} \phi^{p-1} dA
= \Lambda\, H_0(t).
\end {eqnarray} 
Let $l(t)$ be the length of $S(t)$ and let $A(t)$ be the area of 
$\Sigma(t)$. Then combining \eqref{min1} with the isoperimetric inequality 
for minimal surfaces \cite{LSY} and the Cauchy-Schwarz inequality we have 
\[
4\pi A(t) \leq  l^2(t)  \leq  \left ( \int_{S(t)} |\nabla \phi|\, ds \right )
\left ( \int_{S(t)} \frac{ds}{|\nabla \phi|}
\right ) 
=\Lambda\, H_0(t) \int_{S(t)} \frac{ds}{|\nabla \phi|}, 
\]
which we can rearrange to read 
\begin {equation} \label {min2} 
\int_{S(t)} \frac{ds}{|\nabla \phi|} \geq \frac{4\pi A(t)}{\Lambda H_0(t)}.
\end {equation} 
Now, combining \eqref{min0} and \eqref{min2} we have 
\begin {equation} \label {min3}
(H_0^2(t))' = 2 H_0(t) H_0'(t) \leq -\frac{8\pi A(t) t^{p-1}}{\Lambda}. 
\end {equation}

Next we compute 
\begin {eqnarray} \label {min4} 
\frac{dH_1}{dt} & = & -\frac{p}{2}\,\frac{d}{dt} \left( \int_\tau^{\phi_M} 
t^{p-1} \int_{S_\tau} \frac{\langle \nabla \phi, \nabla f \rangle }
{|\nabla \phi|} \,ds \,d\tau\right)  \nonumber \\
& = & -\frac{p\, t^{p-1}}{2} \int_{\Sigma(t)} \Delta_\Sigma f\, dA = 
-p\,t^{p-1} A(t).
\end {eqnarray} 
Combining \eqref{min3} and \eqref{min4} we have 
$$\frac{d}{dt} \left [ H_0^2(t) - \frac{8\pi}{p\Lambda} H_1(t) \right ] \leq 0.$$
We integrate this last inequality from $t=0$ to $t=\phi_M$ and 
use $H_0(\phi_M) = 0 = H_1 (\phi_M)$ to obtain
$$H_0^2 (0) \geq \frac{8\pi}{p\Lambda} H_1 (0).$$
However, the definition of $H_0$ and \eqref{min-H1} tell us that 
$$H_0(0) = \int_\Sigma \phi^{p-1} dA, \qquad 
H_1(0) = \int_\Sigma \phi^p dA,$$
and so, using \eqref{1.2}, we have 
\begin{align*}
\left( \int_\Sigma \phi^{p-1}\,dA\right)^2 
& \geq  \frac{8\pi}{p\Lambda} \int_\Sigma\phi^p\,dA 
= \frac{8\pi}{p} \frac{\left(\int_\Sigma \phi^p\,dA\right)^{(p-2)/p}}{\calCp}\,
\int_\Sigma\phi^p\,dA\\
& = \frac{8\pi}{p\,\calCp} \,\left(\int_\Sigma \phi^p\,dA\right)^{2-2/p}
\end{align*}
which is \eqref{pPR}.

If we have equality in \eqref{pPR} then we must have $4\pi A(t) = 
l^2(t)$ for almost every $t$, which, by Theorem 1 of \cite{LSY} implies 
$\Sigma$ is a flat disk in a two dimensional affine plane. 
\end {proof}

\end{document}